\documentclass[aop,citesort,MSNbibl,dvips]{arximspdf}

\doi{10.1214/10-AOP635}
\volume{40}
\issue{2}
\pubyear{2012}
\firstpage{890}
\lastpage{891}

\begin{document}
\begin{frontmatter}

\title{Correction
Limit theorems for coupled continuous time random walks}\vspace*{12pt}
\runtitle{Correction}
\pdftitle{Correction
Limit theorems for coupled continuous time random walks}

\textit{Ann. Probab.} \textbf{32} (2004) 730--756

\begin{aug}
\author[A]{\fnms{Peter} \snm{Kern}\ead[label=e1]{kern@math.uni-duesseldorf.de}},
\author[B]{\fnms{Mark M.} \snm{Meerschaert}\corref{}\ead[label=e2]{mcubed@stt.msu.edu}}\\ and
\author[C]{\fnms{Hans-Peter} \snm{Scheffler}\ead[label=e3]{scheffler@mathematik.uni-siegen.de}}
\runauthor{P. Kern, M. M. Meerschaert and H.-P. Scheffler}
\affiliation{Heinrich Heine University, Michigan State University
and University~of~Siegen}
\address[A]{P. Kern\\
Mathematisches Institut\\
Heinrich Heine University\\
40225 D\"{u}sseldorf\\
Germany\\
\printead{e1}} 
\address[B]{M. M. Meerschaert\\
Department of Statistics and Probability\\
Michigan State University\\
East Lansing, Michigan 48824\\
USA\\
\printead{e2}}
\address[C]{H.-P. Scheffler\\
Fachbereich Mathematik\\
University of Siegen\\
57068 Siegen\\
Germany\\
\printead{e3}}
\end{aug}

\received{\smonth{8} \syear{2010}}
\revised{\smonth{11} \syear{2010}}



\end{frontmatter}

The converse portion of Theorem 2.2 requires an additional condition,
that the probability measure $\omega$ is such that (2.10) assigns
finite measure to sets bounded away from the origin. The argument on
page 735 must consider~$B_1$ and $B_2$ such that {at least one} is bounded
away from zero, not just the case where both are bounded away from zero.
The condition on~$\omega$ ensures that the integral on page 735 l.--2 is
finite, which is obviously necessary.

The limit process in Theorem 3.4 should read $A(E(t)-)$. If $A(t)$
and~$D(t)$ are dependent, this is a different process than $A(E(t))$. To
clarify the argument, note that
%
\begin{equation}\label{mistake}\qquad
\lim_{h\downarrow0}\frac1hP\{A(s)\in M, s<E(t)\leq s+h\}= P\{A(s-)\in
M | E(t)=s\} p_t(s),
\end{equation}
where $p_t$ is the density of $E(t)$, since $s<E(t)$ in the
conditioning event. For an alternative proof, see Theorem 3.6 in Straka
and Henry \cite{StrakaHenry}. Theorem~4.1 in \cite{coupleCTRW} gives
the density of $A(E(t)-)$.
Examples 5.2--5.6 in \cite{coupleCTRW} provide governing equations for
the CTRW limit process $M(t)=A(E(t)-)$ in some special cases with
simultaneous jumps. Especially, Example 5.5 considers the case where
$Y_i=J_i$ so that $A(t)$ is a stable subordinator and $E(t)=\inf\{
x>0\dvtx A(x)>t\}$ is its inverse or first passage time process. The beta
density for $A(E(t)-)$ given in that example agrees with the result in
Bertoin \cite{bertoin}, page~82. Note that here we have $A(E(t)-)<t$ and
$A(E(t))>t$ almost surely for any $t>0$, by \cite{bertoin}, Chapter III, Theorem~4.


%
\printaddresses

\end{document}